\documentclass[11pt]{article}
\usepackage{amssymb}
\usepackage{amsbsy}
\usepackage{color}
\usepackage{graphicx}
\usepackage{newtxtext}
\usepackage{newtxmath}
\newtheorem{theorem}{Theorem}
\newtheorem{lemma}[theorem]{Lemma}

\newenvironment{proof}{\textbf{Proof.}}{\hfill}

\begin{document}

\title{Stability and eigenvalue bounds for micropolar shear flows}
\author{P. Braz e Silva\thanks{Department of Mathematics, UFPE, Brazil. E-mail: {\tt pablo.braz@ufpe.br}. Partially supported by CAPES, Brazil, \#8881.520205/2020-01 - MATH-AMSUD project 21-MATH-03 (CTMicrAAPDEs), \#88887.311962/2018-00 - PRINT and CNPq, Brazil \#305233/2021-1, \#308758/2018-8, \#432387/2018-8. P. Braz e Silva would like to thank Professor F. Vitoriano e Silva for the 
very fruitful discussions and work on the very first results of this paper.},
J. Carvalho \thanks{Department of Mathematics, UFPE, Brazil. E-mail: {\tt jackellyny.dassy@ufpe.br}. Partially supported CAPES-PRINT, Brazil, \#88887.311962/2018-00}}
\date{}
\maketitle

\begin{abstract}
We prove
eigenvalue bounds for two-dimensional linearized disturbances of parallel flows of micropolar fluids, deriving the Orr-Sommerfeld equations and providing 
a sufficient condition for linear stability of such flows. We also derive wave speed bounds.

\end{abstract}

\section{Introduction}
	In this short note, we study two dimensional perturbations of parallel flows of micropolar, also known as asymmetric, fluids. To this end, we first derive the Orr-Sommerfeld equations for this kind of fluids and, through a variational method, prove some bounds for their eigenvalues. We obtain bounds for the imaginary part of the eigenvalue ensuring stabilty in some regions of the parameters. We also prove some wave speed bounds for perturbations of the base flow.
	 
	 The problem of stability for shear flows of incompressible fluids governed by the Navier-Stokes equations is a classical one. For these equations, Synge \cite{Synge} has stablished eigenvalue bounds assuring linear stability in some region of the parameters $\alpha R$, where $R$ is the Reynolds number and $\alpha$ is the wave number of the perturbations. Later on, Joseph \cite{Joseph} obtained a larger region of stability with respect to $\alpha R$. For the special case of Classical Couette flow, Romanov \cite{Romanov} proved its linear stability for all Reynolds numbers. 
	 
	Here, we are interested in shear flows of micropolar incompressible fluids. Physically, they are fluids consisting of rigid, random oriented (or spherical) particles suspended in a viscous medium, where the deformation of the particles is ignored (for more details see, for example, \cite{Lukas}).
	 Our results here are the generalization to micropolar fluids of the results by Joseph \cite{Joseph} for the classical Navier-Stokes case. For dipolar fluids, similar results have been obtained by Puri \cite{Puri}. 
	
	This work is organized as follows: first, in section 2, following Liu \cite{Liu}, we write the equations for perturbations of parallel flows and derive the Orr-Sommerfeld equations for micropolar fluids. 
	In section 3, we obtain eigenvalue estimates for such equations. 
	The bounds for the imaginary part of the eigenvalues lead to a region of linear stability.
	\section{Linear disturbance equations}
	
	The incompressible micropolar fluids equations are 
	\begin{eqnarray}
	 \mathbf{u}_t + (\mathbf{u} \cdot \nabla ) \mathbf{u} + \nabla p & = & (\mu + \mu_r) \Delta \mathbf{u}+ 2 \mu_r \mbox{curl}\, \mathbf{w}, \label{1'}\\
	 \mathbf{w}_t + (\mathbf{u} \cdot \nabla) \mathbf{w} + 4 \mu_r \mathbf{w} & = &(c_a+c_d) \Delta \mathbf{w} +
	 (c_0+c_d-c_a) \nabla \text{div}\,\mathbf{w} + 2 \mu_r \mbox{curl}\, \mathbf{u}, \label{2'}\\
	 \mbox{div}\, \mathbf{u} & = & 0, \label{3'} 
	\end{eqnarray}
where $\mathbf{u}$ is the linear velocity, $\mathbf{w}$ the angular velocity of particles and $p$ is the pressure. 	
	We are interested in two dimensional perturbations of two dimensional steady parallel flows; so, we look for solutions 
	$\mathbf{u}'$, $\mathbf{v}'$, $p'$ of system \eqref{1'}-\eqref{3'} of the form 
	\begin{eqnarray*}
	\mathbf{u}' & = & \mathbf{U} + \mathbf{u}  \\
	\mathbf{w}' & = & \mathbf{W} + \mathbf{w} \\
	p' & = & P  + p,
	\end{eqnarray*}
where the base flow is $\mathbf{U} = ( U(y) , 0 , 0 )$, $\mathbf{W} = ( 0 , 0 , W(y))$ and the perturbations are 
	$$\textbf{u}(x,y)=(u (x,y),v(x,y),0), \ \ \mathbf{w}(x,y)= (0,0,w (x,y)), $$
	$(x , y ) \in \mathbb{R} \times [0,1]$. 
Componentwise, considering only the first order terms, the equations for the perturbations are 
\begin{eqnarray}
u_t + U u_x + v U' + p_x  & = & (\mu + \mu_r ) \Delta u + 2 \mu_r w_y , \label{5'} \\
 v_t + U v_x + p_y  &= & (\mu + \mu_r ) \Delta v - 2 \mu_r w_x  , \label{6'}\\ 
 w_t + U w_x + v W' + 4 \mu_r w  & = & (c_a + c_d) \Delta w + 2 \mu_r (v_x - u_y), \label{7'}  \\ 
 u_x + v_y  & = & 0  \label{4'}
\end{eqnarray}
with boundary conditions $u = v = w = 0$ at $ y = 0, 1$. 
	Equations \eqref{5'}-\eqref{4'} are the linear disturbance field equations for parallel flow of micropolar fluids. In
	dimensionless form (see \cite{Liu}), these equations read 
	\begin{eqnarray} 
	u_t + U u_{x} + v U'  + p_x & = &  \Bigg(\frac{1}{R_{\mu}} + \frac{1}{2R_{k}}\Bigg)\Delta u  
	+ \frac{R_{0}}{R_{k}} w_{y}, \label{20'} \\
	v_t + U v_{x} + p_{y}& = &   \Bigg(\frac{1}{R_{\mu}} + \frac{1}{2R_{k}}\Bigg)\Delta v 
	- \frac{R_{0}}{R_{k}} w_{x}, \label{21'} \\
	w_t + U w_{x} + v W'   + 2 \frac{R_{0}}{R_{\nu}} w & = & \frac{1}{R_{\gamma}} \Delta w 
	+ \frac{1}{R_{\nu}}(v_{x} - u_{y}), \label{22'} \\
	u_{x} + v_{y} & = & 0, \label{19'}
	\end{eqnarray} 
	where $R_0$, $R_k$, $R_\mu$, $R_\nu$, $R_\gamma$ are dimensionless numbers. 
	Now, 
	consider the dimensionless stream function $\psi$ such that
	\begin{equation} \label{23'}
		u = \psi_{y}, \ \ \ v= - \psi_{x},
	\end{equation}
	and let
	\begin{eqnarray}\label{24'}
		\psi & =&  \varphi(y) e^{i\alpha(x-Ct)}, \\
		 w & = &\omega(y) e^{i\alpha(x-Ct)}, \nonumber
	\end{eqnarray}
	where $\alpha$ is the real wave number and $C$ denotes the mode of disturbances (an equivalent approach is to take the Laplace-Fourier transform of the equations).  
	In order to eliminate $p$ from equations \eqref{20'} and \eqref{21'}, just differentiate \eqref{20'} with respect to $y$, differentiate \eqref{21'} with respect to $x$, subtract one from another and use the incompressibilty condition \eqref{19'} to obtain
	\begin{equation*}
		(\Delta \psi)_{t} + U (\Delta \psi)_{x} - \psi_{x} U'' = \frac{R_{0}}{R_{k}} \Delta w + \Bigg(\frac{1}{R_{\mu}} + \frac{1}{2R_{k}}\Bigg) \Delta^{2} \psi.
	\end{equation*}
Now, using identities \eqref{24'}, we get
	\begin{equation} \label{25'}
		i\alpha[(U - C)(D^{2} - \alpha^{2}) - U''] \varphi = \Bigg(\frac{1}{R_{\mu}} + \frac{1}{2R_{k}}\Bigg) (D^{2} - \alpha^{2})^{2} \varphi + \frac{R_{0}}{R_{k}}(D^{2} - \alpha^{2}) \omega 
	\end{equation}
	where $D= \frac{d}{dy}$.
	Similarly, equation \eqref{22'} becomes 
	\begin{equation*}
		w_{t} + U w_{x} - \psi_{x} W' = -2\frac{R_{0}}{R_{\nu}}w 
		+ \frac{1}{R_{\gamma}}\Delta w 
		- \frac{1}{R_{\nu}}\Delta \psi,
	\end{equation*}
	and, therefore, using identities \eqref{24'}, one obtains
	\begin{equation} \label{26'}
		i\alpha[(U - C)\omega - W' \varphi ] = \Bigg[ \frac{1}{R_{\gamma}}(D^{2} - \alpha^{2}) - \frac{2R_{0}}{R_{\nu}} \Bigg]\omega - \frac{1}{R_{\nu}}(D^{2} - \alpha^{2})\varphi.
	\end{equation}
	The boundary conditions for $\varphi$ and $\omega$ are 
	\begin{equation}\label{bc} 
	\varphi (0)= \varphi(1) = \varphi'(0) = \varphi' (1) = \omega ( 0) = \omega (1) = 0.
	\end{equation}
	Having in mind the classical case, we call equations \eqref{25'} and \eqref{26'}  
	the Orr-Sommerfeld equations for micropolar fluids.
	\section{Eigenvalue bounds}
	
	We are now in position to derive bounds for 
	 the eigenvalues
	$$C=C_{r} +iC_{i}$$
	of the Orr-Sommerfeld problem for micropolar fluids \eqref{25'}, \eqref{26'}, \eqref{bc}. For four times continuously differentiable functions defined on the interval $ y \in [ 0 , 1]$, we denote 
	$\langle f , g \rangle  = \int_0^1 f(y) \overline{g(y)} dy $, $\| f \|^2 = \langle f , f \rangle$ and 
	$\overline{H}$ the complex valued Hilbert space of four times differentiable functions $\varphi$ completed under the norm 
	$\| \varphi'' \|$ by adding limits of sequences satisfying $\varphi ( 0) = \varphi (1) = \varphi' (0) = \varphi' (1) = 0$. 
	Our first result is the following. 
	\begin{lemma}
	The real and imaginary parts of $C=C_{r} +iC_{i}$ satisfy
	\begin{eqnarray} 
		C_{i} & = & \frac{ 
		Q+\overline{Q} 
		- \Big(\frac{1}{\alpha R_{\mu}} + \frac{1}{2\alpha R_{k}}\Big)(||\varphi''||^{2} + 2\alpha^{2}||\varphi'||^{2} + \alpha^{4}||\varphi||^{2})}{||\varphi'||^{2} + \alpha^{2}||\varphi||^{2} + ||\omega||^{2}}  \nonumber \\ 
		& & \mbox{} - \frac{\frac{1}{\alpha R_{\gamma}}(||\omega'||^{2} + \alpha^{2} ||\omega||^{2}) - \frac{2R_{0}}{\alpha R_{\nu}}||\omega||^{2}}{||\varphi'||^{2} + \alpha^{2}||\varphi||^{2} + ||\omega||^{2}} \label{4} \\ 
		& & \mbox{} -\dfrac{\mbox{im}\Big\{ \frac{i R_{0}}{\alpha R_{k}}[\langle\omega,\varphi''\rangle - \alpha^{2} \langle \omega,\varphi \rangle ] + \frac{i}{\alpha R_{\nu}}[\langle \varphi'',\omega \rangle - \alpha^{2} \langle \varphi, \omega \rangle]\Big\}}{||\varphi'||^{2} + \alpha^{2}||\varphi||^{2} + ||\omega||^{2}} , \nonumber
	\end{eqnarray}
	and 
	\begin{equation}\label{5}
		C_{r}  = \dfrac{\int_{0}^{1}[U|\varphi'|^2 + (\alpha^{2} U + \frac{1}{2} U'' )|\varphi |^2 ]\ dy + \int_{0}^{1}U | \omega|^2 \ dy - \frac{1}{2}\int_{0}^{1}W'(\varphi \overline{\omega} + \omega \overline{\varphi}) \ dy}{||\varphi'||^{2} + \alpha^{2}||\varphi||^{2} + ||\omega||^{2}} 
	\end{equation}
	where
	$$Q = \frac{i}{2} \int_{0}^{1}( U' \varphi \overline{\varphi'}  -  W' \varphi \overline{\omega} ) \ dy.$$
\end{lemma}
\begin{proof}
Multipliying equations \eqref{25'} and \eqref{26'} by $\overline{\varphi}$ and $\overline{\omega}$, respectively, and integrating over $[0,1]$, one has 
\begin{eqnarray*}
 & &C ( \| \varphi'\|^2  + \alpha^2  \| \varphi \|^2) +\frac{i}{\alpha} \Big( \frac{1}{R_\mu} + \frac{1}{2R_k} \Big)   ( \| \varphi '' \|^2 + 2 \alpha^2 \| \varphi'\|^2 + \alpha^4 \|\varphi\|^2 )  = \\ 
 & & - \int_0^1 U \varphi'' \overline{\varphi} \, dy +\alpha^2 \int_0^1 U | \varphi |^2 \, dy 
 + \int_0^1 U'' | \varphi |^2 \, dy  - \frac{i R_0}{\alpha R_k} \Big[ \langle w'' , \varphi \rangle - \alpha^2 \langle w , \varphi \rangle \Big] ,
\end{eqnarray*}
and
\begin{eqnarray*}
C \| \omega \|^2 + \frac{i}{\alpha R_\gamma} ( \| \omega ' \|^2 + \alpha^2 \| \omega \|^2 )
+ \frac{2iR_0}{\alpha R_\nu}\|\omega\|^2 & = & \int_0^1 U | \omega |^2 \, dy  - \int_0^1 W' \varphi \overline{\omega} \, dy  \\ & &\mbox{} 
- \frac{i}{\alpha R_\nu} (\langle \varphi'' , \omega \rangle - \alpha^2 \langle \varphi , \omega \rangle ).
\end{eqnarray*}
Adding up these equations and taking the imaginary and real part of the resulting equation, one directly gets the desired indentities \eqref{4} and \eqref{5}.
\end{proof}

From equation \eqref{4}, one has
\begin{eqnarray*} 
		& & \alpha C_{i} (||\varphi'||^{2} + \alpha^{2}||\varphi||^{2} + ||\omega||^{2} ) + \Big(\frac{1}{ R_{\mu}} + \frac{1}{2 R_{k}}\Big)(||\varphi''||^{2} + 2\alpha^{2}||\varphi'||^{2} + \alpha^{4}||\varphi||^{2}) \\ 
& &	\mbox{} + \frac{1}{R_{\gamma}} ||\omega'||^{2}
+ \left(\frac{\alpha^2}{ R_{\gamma}} + \frac{2R_{0}}{ R_{\nu}}\right)||\omega||^{2} = Q + \overline{Q} 
- \left( \frac{ R_{0}}{ R_{k}} + \frac{1}{R_\nu} \right) \left( \mbox{Re} \langle \omega, \varphi'' \rangle -
\alpha^2 \mbox{Re}\langle \varphi , \omega \rangle \right) \\ 
& &
\leq 
 Q + \overline{Q} + \left( \frac{ R_{0}}{ 2R_{k}} + \frac{1}{2R_\nu} \right) (  \| \varphi''\|^2  + \alpha^4 
 \| \varphi \|^2 +2 \| \omega \|^2 ).
	\end{eqnarray*}
	Therefore, if $\frac{1}{ R_{\mu}} + \frac{1}{2 R_{k}} > \frac{ R_{0}}{ 2R_{k}} + \frac{1}{2R_\nu} $ and 
	$\frac{2R_{0}}{ R_{\nu}}> \frac{ R_{0}}{ R_{k}} + \frac{1}{R_\nu}$ (which is assured, for example, by requiring
	$R_1:= \min \{\frac{1}{ R_{\mu}} , \frac{1}{2 R_{k}} ,\frac{R_{0}}{ R_{\nu}} \} > R_2 :=\max
	\{ \frac{ R_{0}}{ 2R_{k}} , \frac{1}{2R_\nu}\}$, one has 
\begin{eqnarray*} 
		& & \alpha C_{i} (||\varphi'||^{2} +  \alpha^{2}||\varphi||^{2} + ||\omega||^{2} ) + 
		(R_1 - R_2 )||\varphi''||^{2} + \Big(\frac{1}{ R_{\mu}} + \frac{1}{2 R_{k}}\Big)2\alpha^{2}||\varphi'||^{2} + (R_1 - R_2)  \alpha^{4}||\varphi||^{2} \\ 
& &	\mbox{} + \frac{1}{R_{\gamma}} ||\omega'||^{2}
+ \left(\frac{\alpha^2}{ R_{\gamma}} + R_1 - R_2 \right) ||\omega||^{2} \leq  Q + \overline{Q} .
	\end{eqnarray*}	
	Simplifying the notation by taking $\frac{1}{R} := \min \{ R_1-R_2 , \frac{1}{ R_{\mu}} + \frac{1}{2 R_{k}},
	\frac{1}{R_\gamma} \} $, we have
	\begin{eqnarray*} 
		& & \hspace{-.5cm}\alpha C_{i} (||\varphi'||^{2} +  \alpha^{2}||\varphi||^{2} + ||\omega||^{2} ) + 
		\frac{1}{R} \left( ||\varphi''||^{2}+ 2\alpha^{2}||\varphi'||^{2} +  \alpha^{4}||\varphi||^{2}  + ||\omega'||^{2}
+ \alpha^2  ||\omega||^{2} + \|\omega \|^2\right) \leq  Q + \overline{Q} ,
	\end{eqnarray*}	
	from which follows the estimate
	\begin{equation} \label{6}
		C_{i} \leq \dfrac{q_{1}\|\varphi'\| \|\varphi\| + q_{2} \|\varphi\| \|\omega\| - (\alpha R)^{-1}(||\varphi''||^{2} + 2\alpha^{2}||\varphi'||^{2} + \alpha^{4}||\varphi||^{2} + ||\omega'||^{2} + \alpha^{2}||\omega||^{2} + ||\omega||^{2})}{||\varphi'||^{2} + \alpha^{2}||\varphi||^{2} + ||\omega||^{2}},
	\end{equation}
	where
	\begin{equation} \label{7}
		q_{1}= \max_{y \in [0,1]}|U'(y)| \mbox{ and } q_{2}= \max_{y \in [0,1]}|W'(y)|.
	\end{equation}
	\begin{theorem}
		Let $R_{1} = \min\left\{\frac{1}{ R_{\mu}} , \frac{1}{2 R_{k}} ,\frac{R_{0}}{ R_{\nu}} \right\}$ and
		$ R_2 =\max \left\{ \frac{ R_{0}}{ 2R_{k}} , \frac{1}{2R_\nu}\right\}$. If $R_{1} > R_{2}$, then 
		\begin{equation} \label{8}
			C_{i} \leq \dfrac{q_{1}+q_{2}}{2\alpha} - \dfrac{\pi^{2} + \alpha^{2}}{\alpha R},
		\end{equation}
		where $\frac{1}{R} = \min \left\{ R_1-R_2 , \frac{1}{ R_{\mu}} + \frac{1}{2 R_{k}},
	\frac{1}{R_\gamma} \right\} $.
		Moreover, no amplified disturbances $(C_{i} > 0)$ of \eqref{25'}, \eqref{26'} and \eqref{bc} exist if 
		\begin{equation} \label{9}
			\alpha R q_{1} < f(\alpha) := max \big\{ M_{1},M_{2}\big\} 
		\end{equation}
		and
		\begin{equation} \label{10}
			\alpha R q_{2} < g(\alpha) := max \big\{ N_{1},N_{2}\big\} ,
		\end{equation}
		where
		\begin{equation} \label{11}
			\left\{
			\begin{array}{rl}
				\displaystyle M_{1} = (4.73)^{2} \pi + 2\alpha^{2} \pi, \\
				\displaystyle M_{2} = (4.73)^{2}\pi + 2^{\frac{3}{2}} \alpha^{3} ,
			\end{array}
			\right.
		\end{equation}
		and
		\begin{equation} \label{12}
			\left\{
			\begin{array}{rl}
				\displaystyle N_{1} = 2(4.73)^{2} \pi + 2\alpha^{3}  , \\
				\displaystyle N_{2} = 2(4.73)^{2}\pi + 2^{\frac{3}{2}} \alpha \pi.
			\end{array}
			\right.
		\end{equation}
	\end{theorem}
	\begin{proof}
		 In the real-valued Hilbert space $H$ corresponding to $\overline{H}$, the following inequalities hold:
		$$||\varphi'||^{2} \geq \pi^{2} ||\varphi||^{2},$$
		\begin{equation} \label{13}
			||\varphi''\|^{2} \geq \pi^{2} ||\varphi'||^{2},
		\end{equation}
		$$||\varphi''\|^{2} \geq (4.73)^{4} ||\varphi||^{2}.$$
	One has 
		\begin{equation} \label{14}
			2 \alpha \|\varphi'\| \|\varphi\| \leq ||\varphi'||^{2} + \alpha^{2}||\varphi||^{2} + ||\omega||^{2},
		\end{equation}
		\begin{equation} \label{15}
			2 \alpha \|\varphi\| \|\omega\| \leq ||\varphi'||^{2} +  \alpha^{2}||\varphi||^{2} + ||\omega||^{2},
		\end{equation}
		\begin{equation} \label{16}
			||\varphi''||^{2} + 2\alpha^{2}||\varphi'||^{2} + \alpha^{4}||\varphi||^{2} + ||\omega'||^{2} + (\alpha^{2} + 1)||\omega||^{2} \geq (\pi^{2} + \alpha^{2})(||\varphi'||^{2} + \alpha^{2}||\varphi||^{2} + ||\omega||^{2}).
		\end{equation}
		Using bounds \eqref{14}, \eqref{15}, \eqref{16} directly in inequality \eqref{6}, one obtains
		\[ 
		C_{i} \leq \frac{q_1 + q_2}{2\alpha} - \frac{\pi^2 + \alpha^2}{\alpha R} .
			\]
			To finish the proof, note that the left hand side of \eqref{6} is non positive if 
			\[
			 \alpha R q_1 \| \varphi' \| \|\varphi \| + \alpha R q_2 \| \varphi \| \|\omega \| \leq ||\varphi''||^{2} + 2\alpha^{2}||\varphi'||^{2} + \alpha^{4}||\varphi||^{2} + ||\omega'||^{2} + \alpha^{2}||\omega||^{2} + ||\omega||^{2}
			\]
			and this surely holds if $\alpha R q_1 \leq \frac{\Gamma_1}{2}$ and $\alpha R q_2 \leq \frac{\Gamma_2}{2}$, where $$\Gamma_{1} := \dfrac{||\varphi''||^{2} + 2\alpha^{2}||\varphi'||^{2} + \alpha^{4}||\varphi||^{2} + ||\omega'||^{2} + (\alpha^{2} + 1)||\omega||^{2}}{\|\varphi\| \|\varphi'\|} $$ and 
			$$\Gamma_{2} := \dfrac{||\varphi''||^{2} + 2\alpha^{2}||\varphi'||^{2} + \alpha^{4}||\varphi||^{2} + ||\omega'||^{2} + (\alpha^{2} + 1)||\omega||^{2}}{\|\varphi\| \|\omega\|}.$$
		Now, just note that
		\begin{eqnarray*}
		\Gamma_{1} & \geq & \frac{\| \varphi''\|^2 + 2 \alpha^2 \| \varphi' \|^2}{\| \varphi '\| \|\varphi\|} \geq (4,73)^2 \pi +  2 \alpha^{2} \pi	= M_1 ,	\end{eqnarray*} 
		\begin{eqnarray*} 
			\Gamma_{1} & \geq & \frac{\| \varphi''\|^2 + 2 \alpha^2 \| \varphi' \|^2 + \alpha^4 \| \varphi \|^2}{\| \varphi '\| \|\varphi\|} \geq (4,73)^2 \pi   \dfrac{2 \alpha^{2}}{\|\varphi\| \|\varphi'\|} \left( \| \varphi'\|^2 +\frac{\alpha^2}{2}\|\varphi\|^2\right)
			\\ & \geq &  (4.73)^{2} \pi + \dfrac{2 \alpha^{2}}{\|\varphi\| \|\varphi'\|} \left( \dfrac{2 \alpha\|\varphi\| \|\varphi'\|}{\surd2}\right)
			\geq (4,73)^2\pi +2^\frac{3}{2} \alpha^3 = M_2 ,\\
		\end{eqnarray*}
		\begin{eqnarray*}
		\Gamma_{2} & \geq & \frac{\| \varphi''\|^2 + \| \omega' \|^2 + \alpha^4 \|\varphi\|^2 +\alpha^2\|\omega\|^2}{\| \varphi \| \omega\|} \geq  \frac{2 \| \varphi''\| \| \omega' \|+ 2\alpha^2 \|\varphi\| \alpha\|\omega\|}{\| \varphi \| \omega\|}\geq
		2 \pi (4,73)^2  +  2 \alpha^{3} 	= N_1 ,	\end{eqnarray*} 
	and
		\begin{eqnarray*} 
			\Gamma_{2} & \geq & 2(4.73)^{2}\pi + 
			\frac{2\alpha^2 \| \varphi'\|^2  + \| \omega \|^2}{\|\varphi\| \|\omega\|}
			\geq 2(4.73)^{2}\pi + 
			\frac{2 \sqrt{2} \alpha \| \varphi'\| \| \omega \|}{\|\varphi\| \|\omega\|} = 2(4.73)^{2}\pi + 
			2^{\frac{3}{2}} \alpha \pi = N_2 ,
			\end{eqnarray*} 
			which proves the result.
			\end{proof}
\begin{figure}
	\includegraphics[scale=0.4]{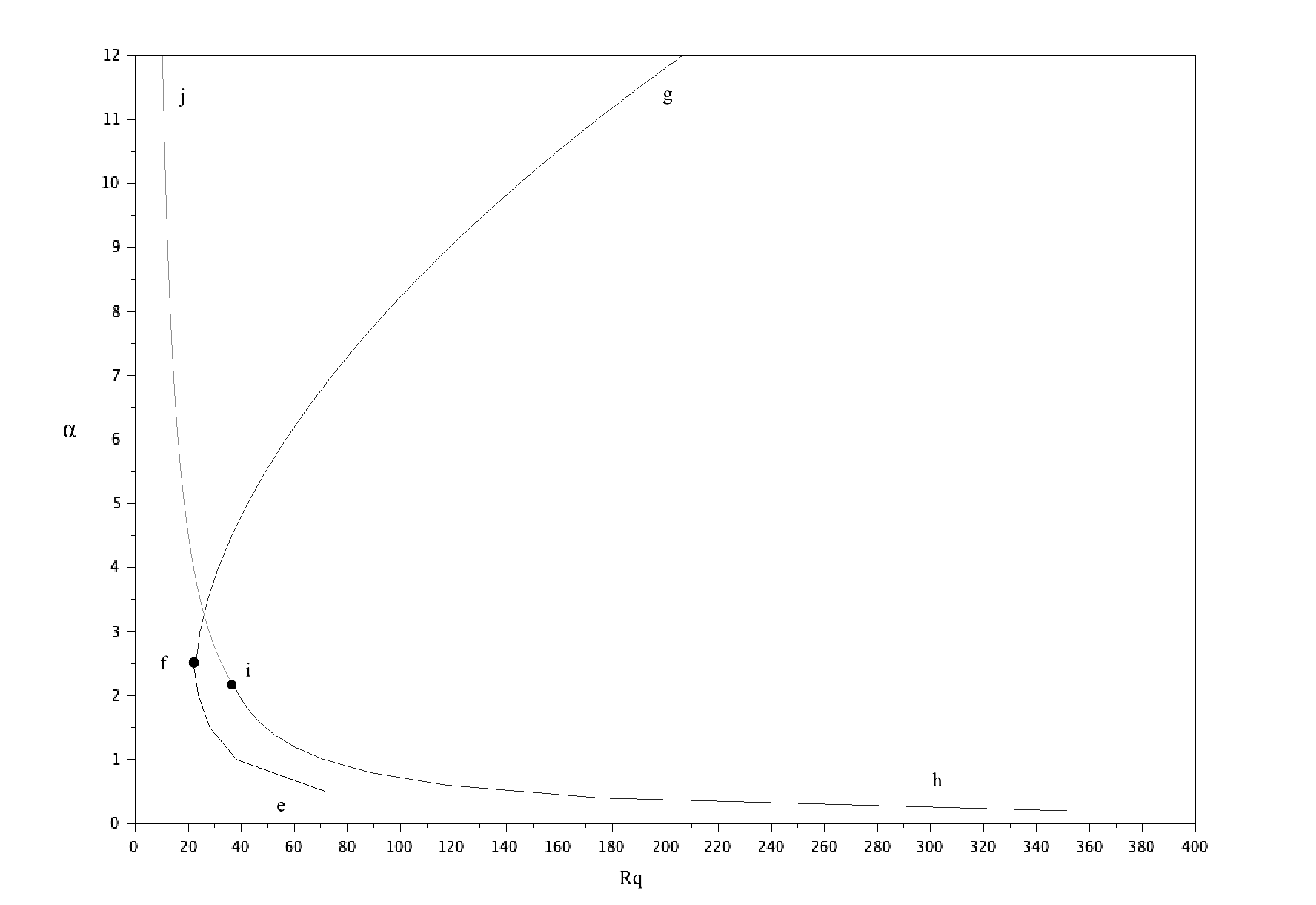}
	\caption{Linear stability bounds for the Orr-Sommerfeld equation for micropolar fluids. The region of certain linear stability lies to the left of the curve \textbf{ej}: \textbf{ef} is the graph of $\frac{M_{2}(\alpha)}{2\alpha}$, \textbf{fg} is the graph of $\frac{M_{1}(\alpha)}{2\alpha}$, \textbf{hi} is the graph of $\frac{N_{2}(\alpha)}{2\alpha}$ and \textbf{ij} is the graph of $\frac{N_{1}(\alpha)}{2\alpha}$.} 
\end{figure}
	\newpage
	Now, we derive estimates for the wave speed $C_{r}$. The result is as follows.
	\begin{theorem} Let $C(\alpha, R) = C_r + i C_i$ be any eigenvalue of \eqref{25'}, \eqref{26'} and \eqref{bc}. Denote 
by $U_{max}, U_{min}, U''_{max}, U''_{min}, W'_{max}$ and $W'_{min}$ the maximum and minimum values on the range of  $U(y), U''(y)$ and $W'(y)$ for $y \in [0,1].$ The following holds:
		\begin{itemize}
		\item[(a)] If $U''_{min} \geq 0 \ \ and \ \ W'_{min} \geq 0$, then
		\begin{equation*} 
			U_{min} - \displaystyle \dfrac{W'_{max}}{2\alpha} \leq C_{r} \leq  U_{max} + \displaystyle \dfrac{U''_{max}}{2(\pi^{2} + \alpha^{2})} + \frac{W'_{max}}{2\alpha}.
		\end{equation*}
		\item[(b)] If $U''_{min} \geq 0 \ \ and \ \ W'_{min} \leq 0 \leq W'_{max}$, then
		\begin{equation*}
			U_{min} - \displaystyle \dfrac{W'_{max}}{\alpha} + \frac{W'_{min}}{2 \alpha} \leq C_{r} \leq U_{max} + \displaystyle \dfrac{U''_{max}}{2(\pi^{2} + \alpha ^{2})} - \frac{W'_{min}}{\alpha} + \frac{W'_{max}}{2\alpha}.
		\end{equation*}
		\item[(c)] If $U''_{min} \geq 0 \ \ and \ \ W'_{max} \leq 0$, then
		\begin{equation*}
			U_{min} + \frac{W'_{min}}{2 \alpha} \leq C_{r} \leq  U_{max} + \displaystyle \dfrac{U''_{max}}{2(\pi^{2} + \alpha ^{2})}  - \frac{W'_{min}}{2 \alpha}.
		\end{equation*}
		\item[(d)] If $ U''_{min} \leq 0 \leq U''_{max} \ \ and \ \ W'_{min} \geq 0$, then
		\begin{equation*}
			U_{min} + \frac{U''_{min}}{2 \alpha^{2}} - \frac{W'_{max}}{2 \alpha} \leq  C_{r} \leq U_{max} + \displaystyle \dfrac{U''_{max}}{2(\pi^{2} + \alpha ^{2})}  + \frac{W'_{max}}{2 \alpha}.
		\end{equation*}
		\item[(e)] If $ U''_{min} \leq 0 \leq U''_{max} \ \ and \ \  W'_{min} \leq 0 \leq W'_{max}$, then
		\begin{equation*}
			U_{min} + \frac{U''_{min}}{2 \alpha^{2}} - \frac{W'_{max}}{\alpha}  + \frac{W'_{min}}{2 \alpha} \leq C_{r} \leq U_{max} + \displaystyle \dfrac{U''_{max}}{2(\pi^{2} + \alpha ^{2})}  - \frac{W'_{min}}{\alpha} + \frac{W'_{max}}{2\alpha}.
		\end{equation*}
		\item[(f)] If $ U''_{min} \leq 0 \leq U''_{max} \ \ and \ \ W'_{max} \leq 0$, then
		\begin{equation*} 
			U_{min} + \displaystyle 
			\dfrac{U''_{min}}{2\alpha^{2}} + \frac{W'_{min}}{2\alpha} \leq C_{r} \leq U_{max} + \displaystyle \dfrac{U''_{max}}{2\alpha^{2}} - \displaystyle \dfrac{W'_{min}}{2\alpha}.
		\end{equation*}
		\item[(g)] If $U''_{max} \leq 0 \ \ and \ \ W'_{min} \geq 0$, then
		\begin{equation*}
			U_{min} + \displaystyle 
			\dfrac{U''_{min}}{2\alpha^{2}}  - \frac{W'_{max}}{2 \alpha} \leq C_{r} \leq U_{max} + \frac{W'_{max}}{2 \alpha}.
		\end{equation*}
		\item[(h)] If $U''_{max} \leq 0 \ \ and \ \ W'_{min} \leq 0 \leq W'_{max}$, then
		\begin{equation*}
			U_{min} + \displaystyle 
			\dfrac{U''_{min}}{2\alpha^{2}} - \frac{W'_{max}}{\alpha}  + \frac{W'_{min}}{2 \alpha}  \leq C_{r} \leq U_{max}  - \frac{W'_{min}}{\alpha} + \frac{W'_{max}}{2\alpha}.
		\end{equation*}
		\item[(i)] If $U''_{max} \leq 0 \ \ and \ \ W'_{max} \leq 0$, then
		\begin{equation*} 
			U_{min} + \displaystyle \dfrac{U''_{min}}{2\alpha^{2}} + \frac{W'_{min}}{2 \alpha}\leq C_{r} \leq U_{max} - \displaystyle \dfrac{W'_{min}}{2\alpha}.
		\end{equation*}
		\end{itemize}
	\end{theorem}
	\begin{proof}
		First of all, note that it follows directly from identity \eqref{5} that
		\begin{equation*}
			C_{r} = U(y_{1}) + \dfrac{ \frac{1}{2}U''(y_{2}) }{ \frac{||\varphi'||^{2}}{||\varphi||^{2}} + \alpha^{2} + \frac{||\omega||^{2}}{||\varphi||^{2}}} - \dfrac{\frac{1}{2}W'(y_{3})[ \langle\varphi, \omega \rangle + \langle \omega,\varphi \rangle ] }{||\varphi'||^{2} + \alpha^{2}||\varphi||^{2} + ||\omega||^{2} }
		\end{equation*}
		where $y_{1}, y_{2}$ and $ y_{3} \in (0,1)$ are mean values. We use this identity and the Poincar\'e inequality 
		\begin{equation} \label{x}
			 ||\varphi'||^{2} \geq \pi^{2} \| \varphi \|^2
		\end{equation}
to prove separately each item of the theorem.
		\begin{enumerate}
			\item[(a)] Case $U''_{min} \geq 0 \ \ and \ \ W'_{min} \geq 0$.
			\begin{eqnarray*}
				C_{r} & = & U(y_{1}) + \dfrac{ \frac{1}{2}U''(y_{2}) }{ \frac{||\varphi'||^{2}}{||\varphi||^{2}} + \alpha^{2} + \frac{||\omega||^{2}}{||\varphi||^{2}}} - \dfrac{\frac{1}{2}W'(y_{3})[ \langle\varphi, \omega \rangle + \langle \omega,\varphi \rangle ] }{||\varphi'||^{2} + \alpha^{2}||\varphi||^{2} + ||\omega||^{2} } \\
				& \leq & U_{max} + \dfrac{ \frac{1}{2}U''_{max} }{ \frac{||\varphi'||^{2}}{||\varphi||^{2}} + \alpha^{2} } - \dfrac{\frac{1}{2 \alpha}W'(y_{3})  [ \langle \alpha \varphi + \omega, \alpha \varphi + \omega \rangle - \alpha^{2}||\varphi||^{2} - ||\omega||^{2}] }{||\varphi'||^{2} + \alpha^{2}||\varphi||^{2} + ||\omega||^{2} } \\
				& \leq & 	U_{max} + \dfrac{ \frac{1}{2}U''_{max} }{ \pi ^{2} + \alpha^{2} } + \dfrac{\frac{1}{2 \alpha}W'(y_{3})[\alpha^{2}||\varphi||^{2} + ||\omega||^{2}]}{||\varphi'||^{2} + \alpha^{2}||\varphi||^{2} + ||\omega||^{2}} \\
				& \leq & U_{max} + \dfrac{ \frac{1}{2}U''_{max} }{ \pi ^{2} + \alpha^{2} } + \dfrac{W'_{max}}{2 \alpha},
			\end{eqnarray*}
			using the Poincar\'e inequality \eqref{x} and that $W'(y_{3})\geq W'_{min}\geq 0$.
			On other hand,
			\begin{eqnarray*}
				C_{r} & = & U(y_{1}) + \dfrac{ \frac{1}{2}U''(y_{2}) }{ \frac{||\varphi'||^{2}}{||\varphi||^{2}} + \alpha^{2} + \frac{||\omega||^{2}}{||\varphi||^{2}}} - \dfrac{\frac{1}{2}W'(y_{3})[ \langle\varphi, \omega \rangle + \langle \omega,\varphi \rangle ] }{||\varphi'||^{2} + \alpha^{2}||\varphi||^{2} + ||\omega||^{2} } \\
				& \geq & U_{min} + \dfrac{ \frac{1}{2}U''_{min} }{ \frac{||\varphi'||^{2}}{||\varphi||^{2}} + \alpha^{2} + \frac{||\omega||^{2}}{||\varphi||^{2}}} 
				- \dfrac{\frac{1}{2}W'(y_{3}) ( \frac{1}{\alpha} \alpha ) 2 \|\varphi\| \|\omega\| }{||\varphi'||^{2} + \alpha^{2}||\varphi||^{2} + ||\omega||^{2} } \\
				& \geq & U_{min} - \dfrac{\frac{1}{2 \alpha}W'(y_{3})(\alpha^{2}||\varphi||^{2} + ||\omega||^{2}) }{||\varphi'||^{2} + \alpha^{2}||\varphi||^{2} + ||\omega||^{2} } \geq U_{min} - \dfrac{W'(y_{3})}{2\alpha}  \geq U_{min} - \dfrac{W'_{max}}{2\alpha}
			\end{eqnarray*}
			using Cauchy's inequality, Young's inequality and the hypothesis $U''_{min} \geq 0.$
			\item[(b)] Case $U''_{min} \geq 0 \ \ and \ \ W'_{min} \leq 0 \leq W'_{max}$.
			\begin{eqnarray*}
				C_{r} & = & U(y_{1}) + \dfrac{ \frac{1}{2}U''(y_{2}) }{ \frac{||\varphi'||^{2}}{||\varphi||^{2}} + \alpha^{2} + \frac{||\omega||^{2}}{||\varphi||^{2}}} - \dfrac{\frac{1}{2}W'(y_{3})[ \langle\varphi, \omega \rangle + \langle \omega,\varphi \rangle ] }{||\varphi'||^{2} + \alpha^{2}||\varphi||^{2} + ||\omega||^{2} } \\
				& \leq & U_{max} + \dfrac{ \frac{1}{2}U''_{max} }{ \frac{||\varphi'||^{2}}{||\varphi||^{2}} + \alpha^{2} } - \dfrac{\frac{1}{2\alpha}W'(y_{3})||\alpha \varphi + \omega||^{2} }{||\varphi'||^{2} + \alpha^{2}||\varphi||^{2} + ||\omega||^{2} } + \dfrac{\frac{1}{2\alpha}W'(y_{3})(\alpha^{2}||\varphi||^{2}+||\omega||^{2})}{||\varphi'||^{2} + \alpha^{2}||\varphi||^{2} + ||\omega||^{2} } \\
				& \leq & U_{max} + \dfrac{ U''_{max} }{ 2(\pi^{2} + \alpha^{2})} - \dfrac{\frac{1}{2\alpha}W'_{min}[||\alpha \varphi + \omega||^{2} + ||\alpha \varphi - \omega||^{2}]}{||\varphi'||^{2} + \alpha^{2}||\varphi||^{2} + ||\omega||^{2} } + \dfrac{W'_{max}}{2 \alpha } \\
				& = & U_{max} + \dfrac{ U''_{max} }{ 2(\pi^{2} + \alpha^{2})} - \dfrac{\frac{1}{\alpha}W'_{min}[\alpha^{2}||\varphi||^{2} + ||\omega||^{2}]}{||\varphi'||^{2} + \alpha^{2}||\varphi||^{2} + ||\omega||^{2} } + \dfrac{W'_{max}}{2 \alpha } \\
				& \leq & U_{max} + \dfrac{ U''_{max} }{ 2(\pi^{2} + \alpha^{2})} - \dfrac{W'_{min}}{\alpha} + \dfrac{W'_{max}}{2 \alpha },
			\end{eqnarray*}
			using \eqref{x}, the hypothesis $W'_{min}\leq 0\leq W'_{max}$ and the parallelogram law. 
			On other hand,
			\begin{eqnarray*}
				C_{r} & = & U(y_{1}) + \dfrac{ \frac{1}{2}U''(y_{2}) }{ \frac{||\varphi'||^{2}}{||\varphi||^{2}} + \alpha^{2} + \frac{||\omega||^{2}}{||\varphi||^{2}}} - \dfrac{\frac{1}{2}W'(y_{3})[ \langle\varphi, \omega \rangle + \langle \omega,\varphi \rangle ] }{||\varphi'||^{2} + \alpha^{2}||\varphi||^{2} + ||\omega||^{2} } \\
				& \geq & U_{min} - \dfrac{\frac{1}{2\alpha}W'_{max}||\alpha \varphi + \omega||^{2} }{||\varphi'||^{2} + \alpha^{2}||\varphi||^{2} + ||\omega||^{2} } + \dfrac{\frac{1}{2\alpha}W'_{min}( \alpha^{2}||\varphi||^{2}+||\omega||^{2})}{||\varphi'||^{2} + \alpha^{2}||\varphi||^{2} + ||\omega||^{2} } \\
				& \geq & U_{min} - \dfrac{\frac{1}{2\alpha}W'_{max}[||\alpha \varphi + \omega||^{2} + ||\alpha \varphi - \omega||^{2}]}{||\varphi'||^{2} + \alpha^{2}||\varphi||^{2} + ||\omega||^{2} } + \dfrac{\frac{1}{2\alpha}W'_{min}(||\varphi'||^{2} + \alpha^{2}||\varphi||^{2}+||\omega||^{2})}{||\varphi'||^{2} + \alpha^{2}||\varphi||^{2} + ||\omega||^{2} } \\
				& \geq & U_{min} - \dfrac{\frac{1}{\alpha}W'_{max}[\alpha^{2}||\varphi||^{2} + ||\omega||^{2}]}{||\varphi'||^{2} + \alpha^{2}||\varphi||^{2} + ||\omega||^{2} } + \dfrac{W'_{min}}{2 \alpha } \\
				& \geq & U_{min} - \dfrac{\frac{1}{\alpha}W'_{max}[||\varphi'||^{2}+\alpha^{2}||\varphi||^{2} + ||\omega||^{2}]}{||\varphi'||^{2} + \alpha^{2}||\varphi||^{2} + ||\omega||^{2} } + \dfrac{W'_{min}}{2 \alpha } \\
				& \geq & U_{min} - \dfrac{W'_{max}}{\alpha} + \dfrac{W'_{min}}{2 \alpha },
			\end{eqnarray*}
			using the hypothesis and the parallelogram law.
			\item[(d)] $ U''_{min} \leq 0 \leq U''_{max} \ \ and \ \ W'_{min} \geq 0$
			\begin{eqnarray*}
				C_{r} & = & U(y_{1}) + \dfrac{ \frac{1}{2}U''(y_{2})||\varphi||^{2} }{ ||\varphi'||^{2} + \alpha^{2} ||\varphi||^{2} + ||\omega||^{2}} - \dfrac{\frac{1}{2}W'(y_{3})[ \langle\varphi, \omega \rangle + \langle \omega,\varphi \rangle ] }{||\varphi'||^{2} + \alpha^{2}||\varphi||^{2} + ||\omega||^{2} } \\
				& > & U_{min}+ \dfrac{ \frac{1}{2}U''_{min} ||\varphi||^{2}}{ ||\varphi'||^{2} + \alpha^{2}||\varphi||^{2} + ||\omega||^{2}   } - \dfrac{\frac{1}{\alpha} \alpha W'(y_{3})\|\varphi\| \|\omega\| }{||\varphi'||^{2} + \alpha^{2}||\varphi||^{2} + ||\omega||^{2} } \\
				& \geq & U_{min}+ \dfrac{ \frac{1}{2}U''_{min} (\frac{1}{\alpha^{2}}||\varphi'||^{2} + ||\varphi||^{2} + \frac{1}{\alpha^{2}}||\omega||^{2})}{ \alpha^{2}( \frac{ ||\varphi'||^{2}}{\alpha^{2}} + ||\varphi||^{2} + \frac{ ||\omega||^{2}}{\alpha^{2}})  }- \dfrac{\frac{1}{2\alpha}W'(y_{3})(\alpha^{2}||\varphi||^{2}+||\omega||^{2})}{||\varphi'||^{2} + \alpha^{2}||\varphi||^{2} + ||\omega||^{2} }  \\
				& \geq & U_{min} + \dfrac{U''_{min}}{2\alpha^{2}} - \dfrac{\frac{1}{2\alpha}W'(y_{3})(||\varphi'||^{2} + \alpha^{2}||\varphi||^{2}+||\omega||^{2})}{||\varphi'||^{2} + \alpha^{2}||\varphi||^{2} + ||\omega||^{2} } \\
				& \geq & U_{min} + \dfrac{U''_{min}}{2\alpha^{2}} - \dfrac{W'_{max}}{2\alpha },
			\end{eqnarray*}
			using that  $U''_{min}\leq 0$, $ W'(y_{3})\geq W'_{min}\geq0$, Cauchy and Young inequalities. The upper bound for $C_r $ follows exactly as in item (a) above. 
		\end{enumerate}
		The proof of the remaining items follow through completely analogous arguments to one of the items proved above.  
		
	\end{proof}


\begin{thebibliography}{99}

%
\bibitem{Joseph} Joseph, D.D., Eigenvalues bounds for the Orr-Sommerfeld equation, {\it J. Fluid Mech.}, 33, pp. 617-621, 1968

\bibitem{Liu} Liu, C.Y., On turbulent flow of micropolar fluids, {\it Internat. J. Engrg. Sci.}, 8, pp. 457-466, 1970
%
\bibitem{Lukas} Lukaszewicz, G.,Micropolar Fluids. Theory and Applications, {\it Modelling and Simulation in Science, Engineering and Technology}, Birkh\"auser, 1999
%
\bibitem{Puri} Puri, P., Stability and eigenvalue bounds of the flow of a dipolar fluid between two parallel plates, {\it Proceedings of the Royal Society A}, 461, pp. 1401-1421, 2005

\bibitem{Romanov} Romanov, V.A., Stability of plane-parallel Couette flow, {\it Functional Analysis and Its Applications}, 7, pp. 137-146, 1973
%
\bibitem{Synge} Synge, J.L., Hydrodynamical stability, {\it Semicentenn. Publ. Amer. Math. Soc.}, 2, pp. 227-269, 1938
%
\end{thebibliography}
\end{document}